\documentclass[english,onecolumn]{IEEEtran}
\usepackage[T1]{fontenc}
\usepackage[latin9]{inputenc}
\usepackage{babel}
\usepackage{textcomp}
\usepackage{amsthm}
\usepackage{amsmath}
\usepackage{amssymb}
\usepackage{graphicx}
\usepackage{esint}
\usepackage[unicode=true]
 {hyperref}
\usepackage{breakurl}

\makeatletter
\usepackage{epstopdf}
\DeclareGraphicsExtensions{.eps}
\providecommand{\tabularnewline}{\\}

  \theoremstyle{plain}
  \newtheorem{lem}{\protect\lemmaname}
  \theoremstyle{plain}
  \newtheorem{thm}{\protect\theoremname}
  \theoremstyle{remark}
  \newtheorem{rem}{\protect\remarkname}

\makeatother

\providecommand{\lemmaname}{Lemma}
\providecommand{\remarkname}{Remark}
\providecommand{\theoremname}{Theorem}

\begin{document}

\title{\textbf{\Large{}Inertial Vector Based Attitude Stabilization of Rigid
Body Without Angular Velocity Measurements}}

\author{L.~Benziane, A.~Benallegue, Y. Chitour, A. Tayebi%
\thanks{This research was partially supported by the iCODE institute, research
project of the Idex Paris-Saclay. L.~Benziane, A. Benallegue are
with LISV, Universit\'{e} de Versailles Saint Quentin, France. Y. Chitour
is with L2S, Universit\'{e} Paris-Sud XI, CNRS and Sup\'{e}lec, Gif-sur-Yvette,
and Team GECO, INRIA Saclay \textendash{} Ile-de-France, France. A.
Tayebi is with the Department of Electrical Engineering, Lakehead
University, Thunder Bay, Ontario, Canada, e-mail: \protect\href{mailto:lotfi.benziane@ens.uvsq.fr, benalleg@lisv.uvsq.fr, yacine.chitour@lss.supelec.fr, atayebi@lakeheadu.ca}{lotfi.benziane@ens.uvsq.fr, benalleg@lisv.uvsq.fr, yacine.chitour@lss.supelec.fr, atayebi@lakeheadu.ca}.%
}}
\maketitle
\begin{abstract}
We address the problem of attitude stabilization of a rigid body,
in which neither the angular velocity nor the instantaneous measurements
of the attitude are used in the feedback, only body vector measurements
are needed. The design of the controller is based on an angular velocity
observer-like system, where a first order linear auxiliary system
based directly on vector measurements is introduced. The introduction
of gain matrices provides more tuning flexibility and better results
compared with existing works. The proposed controller ensures almost
global asymptotic stability. The performance and effectiveness of
the proposed solution are illustrated via simulation results where
the gains of the controller are adjusted using non linear optimization. 
\end{abstract}

\section{Introduction}

Attitude stabilization of rotational motion of rigid body is a classical
problem. Despite the considerable existing solutions, it remains until
today an active research topic. This is due to the large field of
applications such as robotics, unmanned aerial vehicles (UAVs), satellites,
marine vehicles, etc. The problem of attitude control was treated
using several type of parametrization of the attitude \cite{Shuster1993}.
Classical solutions to this problem have been proposed using local-minimal
parametrization which lies in $\mathbb{R}^{3}$, such as Euler angles
or modified Rodriguez parameters (see, for instance, \cite{Pflimlina2010,Zou2010,Castillo2006,Wen1991}).
The global-unique representation, which is the natural parametrization
of the attitude is the direction cosine rotation matrix that lies
in the special orthogonal group $SO(3)$. As a consequence, many recent
solutions use this parametrization (see, for instance, \cite{Bayadi2014,Lee2013,Chaturvedi2011,Mahony2008}).
However, for simplicity of analysis and numerical implementation reasons,
a considerable number of solutions to the problem of attitude stabilization
of rigid bodies rather use quaternion parametrization as global representation
which lies in the unit sphere $\mathbb{S}^{3}$ (see, for instance,
\cite{Fresk2013,Zhu2011,Mayhew2009,Tayebi2006,Joshi1995}).

Since attitude control and stabilization is as an interesting theoretical
and technical problem, many scenarii were studied in the literature
(see for instance \cite{Nordkvist2010,Pounds2007}, \cite{Tayebi2007},\cite{Morin1997,Krishnan1994}
and \cite{Byrnes1991}). The interesting and challenging scenario
is the attitude stabilization without angular velocity. In fact, the
main goal is to stabilize the attitude without the use of gyroscopes,
which can be very expensive or vital to the system, like gyroscopes
on Hubble used for pointing the telescope. They measure attitude when
Hubble is changing its pointing from one target (a star or planet)
to another, and they help control the telescope's pointing while scientists
are observing targets. There are a total of six gyroscopes on board--three
serve as backups. In 2009, all six of Hubble's gyroscopes had to be
replaced and one can imagine the cost generated. At the light of these
problems, it is conceivable to reduce costs and ensure continuity
of the mission of the rigid body despite the failure of the gyroscopes
when this type of controllers is used. Many works in the literature
dealt with attitude stabilization without angular velocity problem
(see, for instance \cite{Thakur2014,Benziane2014,Filipe2013,Tayebi2013,Xiao2013,Schlanbusch2012,Tayebi2007a,Akella2001,Wong2001}),
some of them exploited the passivity of the system such as \cite{Tayebi2008,Costic2001,Tsiotras1998,Lizarralde1996,Egeland1994}.

In almost all results dealing with the case of attitude control without
angular velocity, the ``instantaneous measurements of the attitude''
are used in the control law. As there is no sensor which directly
measures the attitude of a rigid body, the aforementioned velocity-free
controllers require some kind of attitude observer relying on the
available direction sensors. However, all static algorithms based
only on body vector measurements are very sensitive to noise. Also,
all the most efficient dynamic attitude estimation algorithms make
use of the body measurements and the angular velocity information
to estimate the attitude of the rigid body. To overcome this problem,
a velocity-free attitude control scheme, that incorporates explicitly
vector measurements instead of the attitude itself, has been proposed
for the first time in \cite{Tayebi2013}. As claimed in \cite{Tayebi2013},
this class of controllers can be qualified as the class of \textit{true
velocity-free attitude controller}s.

Since it is impossible to achieve a global asymptotic stabilization
using continuous time invariant state feedback \cite{Bhat2000}, the
attitude control scheme presented in this work makes use the notion
of ``Almost Global Asymptotic Stabilization'' (AGAS) of the closed
loop system. Therefore, this work and that proposed in \cite{Tayebi2013}
present a stronger stability property compared to \cite{Thakur2014},
where the convergence depends on a non trivial condition on initial
conditions.

The proposed solution given in this paper can be regarded as an continuation
of \cite{Tayebi2013}. The main differences are the following: (a)
the use of an auxiliary system in terms of body vector measurements,
defined on $\mathbb{R}^{3}$, rather than that of an auxiliary system
defined on $\mathbb{S}^{3}$; (b) the explicit design of an angular
velocity observer which is used in the design of the stabilizing feedback.
As a consequence, the set of unstable equilibria of the closed loop
dynamics of our auxiliary error system is reduced as compared to that
of \cite{Tayebi2013}. It is also shown that our auxiliary error system
does not make use the inertial fixed reference vectors as in \cite{Tayebi2013}.
The quaternion parametrization is used in the main analysis and the
final results are rewritten with rotations expressed in $SO(3)$ by
simple projection. We also show that the introduction of gain matrices
improves drastically the controller performance with respect to both
\cite{Tayebi2013} and \cite{Thakur2014}. Moreover and contrarily
to what is stated in \cite{Bayadi2014,Tayebi2013} we prove that the
set of control gains leading to a continuum of equilibria of the closed
loop system is an algebraic variety of positive co-dimension, independently
on the choice of the observed vectors. Finally, in order to adjust
properly the controller gains, we rely a non-linear optimal tuning
method.

The result presented in this paper extends those from \cite{Benziane2014}
where a scalar gain was used in the control law. In addition, a complete
and rigorous mathematical analysis is presented in this version.

\section{Notations and Problem Formulation\label{sec:Preliminary,-assumptions-and}}

\subsection{Notations}

To perform a rotation in Euclidean space, we used either a rotation
matrix $R$ or a unit-quaternion $Q=[q_{0},q^{T}]^{T}$. We assume
that $R\in SO(3)=\{R\in\mathcal{\mathrm{\mathbb{R}}}^{3\times3}\mid R^{T}R=RR^{T}=I_{3\times3},det(R)=1\}$
and $Q\in\mathbb{S}^{3}=\{Q\in\mathcal{\mathrm{\mathbb{R}}}^{4}\mid Q^{T}Q=1\}$.
The multiplication of two quaternions $P=(p_{0},p^{T})^{T}$ and $Q=(q_{0},q^{T})^{T}$
is denoted by ``$\odot$'' and defined as $P\odot Q=\left[\begin{array}{c}
p_{0}q_{0}-p^{T}q\\
p_{0}q+q_{0}p+p\times q
\end{array}\right].$ We use $\mathfrak{so}(3)$ to denote the Lie algebra of $SO(3)$,
i.e., the set of skew symmetric matrices and we set $S$ as the Lie
algebra isomorphism from $\mathbb{R}^{3}\rightarrow\mathfrak{so}(3)$
which associates to $x=[x_{1},x_{2},x_{3}]^{T}$ the skew-symmetric
matrix $S(x)$ given by 
\[
S(x)=\left[\begin{array}{ccc}
0 & -x_{3} & x_{2}\\
x_{3} & 0 & -x_{1}\\
-x_{2} & x_{1} & 0
\end{array}\right].
\]
Note that for every $x,y\in\mathbb{R}^{3}$, one has $S(x)y=x\times y$
where $\times$ stands for the vector cross product.

The mapping $\mathcal{R}:\mathbb{S}^{3}\rightarrow SO(3)$ given by
Rodrigues' rotation formula \cite{Hestenes1999}

\begin{equation}
\begin{array}{c}
\mathcal{R}(Q)=I_{3\times3}+2q_{0}S(q)+2S(q)^{2},\end{array}\label{eq: rodriguez}
\end{equation}
defines a double covering map of $SO(3)$ by $\mathbb{S}^{3}$, i.e.,
for every $R\in SO(3)$ the equation $\mathcal{R}(Q)=R$ admits exactly
two solutions $Q_{R}$ and $-Q_{R}$. As a consequence, a vector field
$f$ of $\mathbb{S}^{3}$ projects onto a vector field of $SO(3)$
if and only if, for every $Q\in\mathbb{S}^{3}$, $f(-Q)=-f(Q)$ (where
we have made the obvious identification between $T_{Q}\mathbb{S}^{3}$
the tangent space of $\mathbb{S}^{3}$ at $Q$ and $T_{-Q}\mathbb{S}^{3}$
the tangent space of $\mathbb{S}^{3}$ at $-Q$).

In what follows and for simplicity, the notations below are used. 
\begin{itemize}
\item If $m$ is a positive integer, $M_{m}(\mathbb{R})$ is used to denote
the set of $m$ by $m$ matrices with real entries; $\boldsymbol{0}_{3}$,
$\boldsymbol{0}_{3n}$, $\boldsymbol{0}$ and $\boldsymbol{I}$ denote
the $3$ by $3$ zero matrix, the $3n$ by $1$ zero vector, the $3$
by $1$ zero vector and the $3$ by $3$ identity matrix respectively; 
\item $\left\{ \mathcal{B}\right\} $ and $\left\{ \mathcal{I}\right\} $
denote an orthonormal body-attached frame with its origin at the center
of gravity of the rigid-body and the inertial reference frame on earth
respectively. 
\end{itemize}
For every $x,y\in\mathbb{R}^{3}$ and a given $R\in SO(3)$ one has
the following \cite{Markley2014} 
\[
\begin{array}{rrl}
S(x)y & = & -S(y)x,\ S(x)x=\boldsymbol{0},\\
S(x)S(y) & = & yx^{T}-x^{T}y\boldsymbol{I},\ S^{2}(x)=xx^{T}-x^{T}x\boldsymbol{I},\\
S(S(x)y) & = & S(x)S(y)-S(y)S(x),\ S(Rx)=RS(x)R^{T}.
\end{array}
\]

\subsection{Problem formulation}

Let $n\geq2$ be an integer. The attitude kinematics of a rigid body
in 3D space is given by

\begin{equation}
\dot{R}(t)=R(t)S(\omega(t)),\label{eq:Attitude kinematics R}
\end{equation}
where $R\in SO(3)$. The equivalent kinematics evolving in $\mathbb{S}^{3}$
are given by 
\begin{equation}
\dot{Q}(t)=\left[\begin{array}{c}
\dot{q}_{0}(t)\\
\dot{q}(t)
\end{array}\right]=\left[\begin{array}{c}
-\frac{1}{2}q^{T}(t)\omega(t)\\
\frac{1}{2}(q_{0}(t)\boldsymbol{I}+S(q(t)))\omega(t)
\end{array}\right],\label{eq:real kinematics}
\end{equation}
where $\omega(t)$ being the angular velocity of the rigid body expressed
in $\left\{ \mathcal{B}\right\} $ and $Q\in\mathbb{S}^{3}$ is the
unit quaternion. Let $b_{i}(Q(t))\in\mathcal{\mathrm{\mathbb{R}}}^{3}$
$(i=1,\cdots,n)$ be a measured vector expressed in $\left\{ \mathcal{B}\right\} $.
The relation between $b_{i}(t)$ and its corresponding fixed inertial
vector $r_{i}\in\mathcal{\mathrm{\mathbb{R}}}^{3}$ are given by

\begin{equation}
b_{i}(Q(t))=R^{T}(t)r_{i}\label{eq:relation bi ri}
\end{equation}
As a consequence we have $b_{i}(-Q)=b_{i}(Q)$ for $1\leq i\leq n$
and $Q\in\mathbb{S}^{3}$.

Using (\ref{eq:Attitude kinematics R}) and (\ref{eq:relation bi ri}),
one can get the reduced attitude kinematics

\begin{equation}
\dot{b}_{i}(Q(t))=-S(\omega(t))b_{i}(Q(t)),\, i=1,\cdots,n.\label{eq:dynamic de bi}
\end{equation}

The simplified rigid body rotational dynamics are governed by 
\begin{equation}
\begin{array}{c}
J\dot{\omega}(t)=-S(\omega(t))J\omega(t)+\tau(t),\end{array}\label{eq. real dynamics}
\end{equation}
where 
\begin{itemize}
\item $J\in\mathbb{R}^{3\times3}$ is a symmetric positive definite constant
inertia matrix of the rigid body expressed in $\left\{ \mathcal{B}\right\} $; 
\item $\tau(t)$ is the external torque applied to the system expressed
in $\left\{ \mathcal{B}\right\} $; 
\item $\omega(t)$ being the angular velocity of the rigid body expressed
in $\left\{ \mathcal{B}\right\} $. 
\end{itemize}
The problem addressed in this work is the design of an attitude stabilization
control $\tau(t)$ based only on inertial measurements $b_{i}(t)$,
\textit{without using the angular velocity $\omega(t)$ in the feedback}.

\subsection{Assumptions\label{sub:Assumtions}}

We make the following assumptions for the rest of the paper. 
\begin{description}
\item [{$A1$}] We assume that only the $n$ vector-valued functions of
time $b_{i}(t)$ are measured and we do not make any similar assumption
on angular velocity vector $\omega(t)$. Moreover, note that the $b_{i}$'s
actually depend on the rotation $R$ and one could also write them
as $b_{i}(R(t))$ or $b_{i}(Q(t))$ if we choose quaternions instead
of rotations. In the sequel, we will write either $b_{i}(t)$ or $b_{i}(Q(t))$. 
\item [{$A2$}] At least two vectors $r_{1}$, $r_{2}$ are non collinear.
As a consequence, $b_{1}(t)$ and $b_{2}(t)$ are linearly independent
for all non negative times. 
\item [{$A3$}] The desired rigid body attitude is defined by the constant
rotation matrix $R_{d}$, relates an inertial vector $r_{i}$ to its
corresponding vector in the desired frame, \textit{i.e.}, $b_{i}^{d}=R_{d}^{T}r_{i}$,
with $\dot{b}_{i}^{d}=\boldsymbol{0}$. An equivalent constant desired
unit-quaternion $Q_{d}$ is defined as $R_{d}=\mathcal{R}(Q_{d})$. 
\end{description}

\section{Handling the lack of angular velocity and Design of the attitude
controller\label{sec:Design-of-the}}

\subsection{Angular velocity observer-like system\label{sub:Handling-the-lack}}

As well known, the reduced attitude kinematic is defined by (\ref{eq:dynamic de bi}).
Define $\Gamma=diag(\Lambda_{1},\cdots,\Lambda_{n})$, where $\Lambda_{i}$
is a symmetric positive definite $3\times3$ matrix, for $1\leq i\leq n$.
Define the symmetric matrix $M(t)=\sum_{i=1}^{n}S(b_{i}(t))^{T}\Lambda_{i}S(b_{i}(t))$,
which is positive definite thanks to Assumption $A2$.

Multiplying (\ref{eq:dynamic de bi}) by $S(b_{i}(t))\Lambda_{i}$
for $1\leq i\leq n$ and doing the sum gives 
\begin{equation}
\sum_{i=1}^{n}S(b_{i}(t))\Lambda_{i}\dot{b}_{i}(t)=-M(t)\omega(t).\label{eq: omega_real first}
\end{equation}
From (\ref{eq: omega_real first}) the true angular velocity $\omega(t)$
is given by 
\begin{equation}
\omega(t)=-M^{-1}(t)\sum_{i=1}^{n}S(b_{i}(t))\Lambda_{i}\dot{b}_{i(t)}\label{eq:true angular velocity}
\end{equation}

Since $\dot{b}_{i}(t)$ is not a measured quantity, we propose the
following new angular velocity observer-like signal

\begin{flushleft}
\begin{equation}
\hat{\omega}(t)=-M^{-1}(t)\sum_{i=1}^{n}S(b_{i}(t))\Lambda_{i}\dot{\hat{b}}_{i}(t),\label{eq: Angular velocity observer}
\end{equation}
where the vector $\dot{\hat{b}}_{i}(t)$ can be viewed as an estimate
of the vector $\dot{b}_{i}(t)$ using the following linear first-order
filter on $b_{i}$ ($i=1,\cdots,n$). 
\par\end{flushleft}

\begin{flushleft}
\begin{equation}
\dot{\hat{b}}_{i}(t)=A_{i}(b_{i}(t)-\hat{b}_{i}(t)),\label{eq: dynamic estimator}
\end{equation}
where the constant matrices $A_{i}\in\mathbb{R}^{3\times3}$ are chosen
as $A_{i}=R_{d}^{T}P_{i}(\Lambda_{i})R_{d}$, for $1\leq i\leq n$,
with $P_{i}$ polynomial of degree two which is positive on $\mathbb{R}_{+}^{\ast}$.
As a trivial consequence, one deduces that, for $1\leq i\leq n$,
$R_{d}A_{i}R_{d}^{T}$ is symmetric positive definite and commutes
with $\Lambda_{i}$. Set $A=diag(A_{1},\cdots,A_{n})$ and $A_{d}=diag(R_{d}A_{1}R_{d}^{T},\cdots,R_{d}A_{n}R_{d}^{T})$.
Then $\Gamma$ and $A_{d}$ commute. 
\par\end{flushleft}

We define an error for the linear first-order filter by $\tilde{b}_{i}(t)=b_{i}(t)-\hat{b}_{i}(t)$.
Using (\ref{eq: dynamic estimator}), (\ref{eq:dynamic de bi}) leads
to the following error dynamics $\dot{\tilde{b}}_{i}(t)=-A_{i}\tilde{b}_{i}(t)+S(b_{i}(t))\omega,$
which can be rewritten using the state vector defined by $\xi(t):=[\tilde{b}_{1}^{T}(t),\cdots,\,\tilde{b}_{n}^{T}(t)]^{T}$,
as 
\begin{equation}
\dot{\xi}(t)=-A\xi(t)+B(t)\omega(t),\label{eq:new auxiliary error dynamics}
\end{equation}
where $B(t)=\left[\begin{array}{ccc}
S(b_{1}(t))^{T} & \cdots & S(b_{n}(t))^{T}\end{array}\right]^{T}$.

Finally, the angular velocity observer-like signal can be written
as

\begin{equation}
\hat{\omega}(t)=M^{-1}(t)B^{T}(t)\Gamma A\xi(t).\label{eq:omega hat observer}
\end{equation}

\subsection{Controller Design\label{sub:Controller-Design}}

First, the orientation error is defined by

\begin{equation}
\bar{R}(t)=R(t)R_{d}^{T},\label{eq:Error def in R}
\end{equation}
where $R(t)$ is a rotation matrix and $R_{d}$ is a constant desired
rotation matrix. From (\ref{eq:Attitude kinematics R}) and (\ref{eq:Error def in R})
one can obtain the attitude dynamics errors in term of matrix rotation
as follows

\begin{equation}
\dot{\bar{R}}(t)=\bar{R}(t)S(R_{d}\omega(t)),\label{eq:Error dynamics in SO3-1}
\end{equation}
$\bar{R}(t)$ correspond to the quaternion error $\bar{Q}(t)=Q(t)\odot Q_{d}^{-1}(t)\equiv[\bar{q}_{0}(t),\bar{q}(t)^{T}]^{T}$
whose dynamics is governed by 
\begin{equation}
\left[\begin{array}{c}
\dot{\bar{q}}_{0}(t)\\
\dot{\bar{q}}(t)
\end{array}\right]=\left[\begin{array}{c}
-\frac{1}{2}\bar{q}^{T}(t)R_{d}\omega(t)\\
\frac{1}{2}\left(\bar{q}_{0}(t)\boldsymbol{I}+S\left(\bar{q}(t)\right)\right)R_{d}\omega(t)
\end{array}\right],\label{eq:erreur orientation quat}
\end{equation}

The reduced orientation error is given by $\bar{b}_{i}(\bar{Q}(t))=b_{i}(Q(t))-b_{i}^{d}$.
Therefore, on can get 
\begin{equation}
\bar{b}_{i}(\bar{Q}(t))=R_{d}^{T}(\bar{R}(t)^{T}-I)r_{i},\label{eq:bi_bar_dynamics SO(3)}
\end{equation}
where $1\leq i\leq n$ which can be rewritten using (\ref{eq: rodriguez})
as

\begin{equation}
\bar{b}_{i}(\bar{Q}(t))=-2R_{d}^{T}(\bar{q}_{0}(t)\boldsymbol{I}-S(\bar{q}(t)))S(\bar{q}(t))r_{i}.\label{eq:bi_bar quaternion}
\end{equation}

We propose the following control law 
\begin{equation}
\tau(t)=z_{\rho}(t)-M\hat{\omega}(t),\label{eq:control law}
\end{equation}
where the term $z_{\rho}(\cdot)$ was introduced in \cite{Tayebi2013}
and is given by 
\begin{equation}
z_{\rho}(t)=\sum_{i=1}^{n}\rho_{i}S(b_{i}^{d})b_{i},\label{eq11}
\end{equation}
where the coefficients $\rho_{i}$'s are arbitrary positive constants.
Define 
\begin{equation}
W_{\rho}=-\sum_{i=1}^{n}\rho_{i}S^{2}(r_{i}),\label{eq:Wrho def}
\end{equation}
The matrix $W_{\rho}$ is positive definite, see Lemma 2 of \cite{Tayebi2013}.
Then, it has been shown in Lemma 1 of \cite{Tayebi2013} that one
can actually rewrite $z_{\rho}(\cdot)$ as

\begin{equation}
z_{\rho}(t)=-2R_{d}^{T}(\bar{q}_{0}(t)\boldsymbol{I}-S(\bar{q}(t)))W_{\rho}\bar{q}(t).\label{eq:Expression of z_rho}
\end{equation}
One finally gets that the controller $\tau(\cdot)$ can be expressed
as

\begin{equation}
\tau(t)=-2R_{d}^{T}(\bar{q}_{0}(t)\boldsymbol{I}-S(\bar{q}(t)))W_{\rho}\bar{q}(t)-M\hat{\omega}(t).\label{eq:  control law quaternion}
\end{equation}

Using (\ref{eq:new auxiliary error dynamics}), (\ref{eq:erreur orientation quat}),
(\ref{eq. real dynamics}) and (\ref{eq:  control law quaternion}),
we obtain the following closed loop dynamics 
\begin{equation}
\begin{cases}
\begin{array}{lcl}
\dot{\xi} & = & -A\xi+B(\bar{Q})\omega,\\
\dot{\bar{q}}_{0} & = & -\frac{1}{2}\bar{q}^{T}R_{d}\omega,\\
\dot{\bar{q}} & = & \frac{1}{2}(\bar{q}_{0}\boldsymbol{I}+S(\bar{q}))R_{d}\omega,\\
J\dot{\omega} & = & -S(\omega)J\omega-2R_{d}^{T}(\bar{q}_{0}\boldsymbol{I}-S(\bar{q}))W_{\rho}\bar{q}-M\hat{\omega}.
\end{array}\end{cases},\label{eq: closed loop dynamics quaternion0}
\end{equation}
One can make a further simplification by changing variables as follows:
\[
\xi\ \rightarrow[R_{d}\tilde{b}_{1}^{T}(Q(t)),\cdots,R_{d}\tilde{b}_{n}^{T}(Q(t))]^{T},\quad\omega\ \rightarrow R_{d}\omega.
\]
By setting 
\[
J_{d}:=R_{d}JR_{d}^{T},\ B_{d}:=\left[\begin{array}{ccc}
S(R_{d}b_{1})^{T} & \cdots & S(R_{d}b_{n})^{T}\end{array}\right]^{T},
\]
and by making obvious abuse of notations (i.e., we keep the variables
$\xi$ and $\omega$) we end up with an autonomous differential equation
\begin{equation}
\begin{cases}
\begin{array}{lcl}
\dot{\xi} & = & -A_{d}\xi+B_{d}\omega,\\
\dot{\bar{q}}_{0} & = & -\frac{1}{2}\bar{q}^{T}\omega,\\
\dot{\bar{q}} & = & \frac{1}{2}(\bar{q}_{0}\boldsymbol{I}+S(\bar{q}))\omega,\\
J_{d}\dot{\omega} & = & -S(\omega)J_{d}\omega-2(\bar{q}_{0}\boldsymbol{I}-S(\bar{q}))W_{\rho}\bar{q}-B_{d}^{T}\Gamma A_{d}\xi.
\end{array}\end{cases},\label{eq: closed loop dynamics quaternion}
\end{equation}
Note that $J_{d}$ is a real symmetric positive definite matrix. If
one defines the state $\chi:=(\xi,\,\bar{Q},\,\omega)$ where $\bar{Q}\equiv\left[\begin{array}{c}
\bar{q}_{0}\\
\bar{q}
\end{array}\right]\in\mathbb{S}^{3}$ and the state space $\Upsilon:=\mathbb{R}^{3n}\times\mathbb{S}^{3}\times\mathbb{R}^{3}$,
one can rewrite (\ref{eq: closed loop dynamics quaternion}) as $\dot{\chi}=F(\chi)$
where $F$ gathers the right-hand side of (\ref{eq: closed loop dynamics quaternion})
and defines a smooth vector field on $\Upsilon$. Moreover, note that
$\bar{Q}$ and $-\bar{Q}$ represents the same physical rotation,
implying that (\ref{eq: closed loop dynamics quaternion}) projects
on $SO(3)$ as an autonomous differential equation. We will use that
fact in Subsection~\ref{sub:Almost-Global-Asymptotic}. 
\begin{lem}
\label{lem:generic0} With the notations above, one gets that the
matrix $W_{\rho}$ defined in \eqref{eq:Wrho def} has simple eigenvalues
generically with respect to $\rho=(\rho_{1},\cdots,\rho_{n})\in(\mathbb{R}_{+}^{\ast})^{n}$. \end{lem}
\begin{IEEEproof}
For $\rho\in(\mathbb{R}_{+}^{\ast})^{n}$, let $P_{\rho}(\cdot)$
be the characteristic polynomial of $W_{\rho}$ and $\Delta(\rho)$
its discriminant \cite{Smith1998}. Recall that $\Delta(\rho)=0$
if and only if $P_{\rho}(\cdot)$ admits a multiple root. Since $W_{\rho}$
is a $3$ by $3$ real symmetric positive definite matrix or every
$\rho\in(\mathbb{R}_{+}^{\ast})^{n}$, $\Delta(\rho)$ is actually
a homogeneous polynomial of degree four in $\rho$. Thus the locus
$\Delta(\rho)=0$ defines an algebraic variety of co-dimension one
in $(\mathbb{R}_{+}^{\ast})^{n}$ and, on its complementary set $\mathcal{S}$
in $(\mathbb{R}_{+}^{\ast})^{n}$, $W_{\rho}$ has simple eigenvalues. 
\end{IEEEproof}
This genericity result serves a justification to the following working
hypothesis, which will hold for the rest of the paper.

\[
\hbox{\textbf{(GEN)}\;\ensuremath{W_{\rho}}has\; simple\; eigenvalues.}
\]

\section{Stability Analysis of the Proposed Controller\label{sec:Stability-Analysis}}

In this section, we present a rigorous analysis using two representations.
As often, it turns out that it is simpler for the stability analysis
to use unit quaternions for the representation of rotations instead
of elements of $SO(3)$, even-though we are ultimately interested
in a result formulated in terms of orthogonal matrices. This is why
we first complete the stability analysis and obtain a first theorem
(Theorem~\ref{thm: theorem RSR}) using unit quaternions and, in
a second step, we state our main result in terms of of elements of
$SO(3)$ by simply projecting Theorem~\ref{thm: theorem RSR} using
Rodriguez formula \eqref{eq: rodriguez}. 
\begin{lem}
\label{lem:The-solutions-of z_rho}Under the hypothesis \textbf{(GEN)},
the solutions of equation $z_{\rho}=\boldsymbol{0}$ where \textup{$z_{\rho}$
}is defined by (\ref{eq:Expression of z_rho}) are the following:
$(a)$ the two points $\pm(1,\boldsymbol{0})$; the six points $\pm(0,v_{i})$,
$1\leq i\leq3$, with $(v_{1},v_{2},v_{3})$ being an orthonormal
basis diagonalizing $W_{\rho}$. \end{lem}
\begin{IEEEproof}
Let $(q_{0},q)\in\mathbb{S}^{3}$ such that $z_{\rho}=\boldsymbol{0}$,
i.e., 
\[
(q_{0}\boldsymbol{I}-S(q))W_{\rho}q=0.
\]
If $q_{0}\neq0$, it is immediate to see that $q_{0}\boldsymbol{I}-S(q)$
is invertible and thus $q=\boldsymbol{0}$, finally implying that
$q_{0}=\pm1$. If $q_{0}=0$, we are left with the equation $S(q)W_{\rho}q=\boldsymbol{0}$.
According to the properties of $S(q)$ with $q\in\mathbb{S}^{2}$,
we get that $q$ is an eigenvector of $W_{\rho}$ with unit length.
We conclude by using \textbf{(GEN)}. 
\end{IEEEproof}
Consider the following non negative differentiable function $V:\:\Upsilon\rightarrow\mathbb{R}^{+}$
\begin{equation}
V=\xi^{T}\Gamma A_{d}\xi+4\bar{q}^{T}W_{\rho}\bar{q}+\omega^{T}J_{d}\omega,\label{eq: Lyapunov function}
\end{equation}
which is radially unbounded over $\Upsilon$ since $W_{\rho}$ and
$J_{d}$ are positive definite. Moreover, since $\Gamma$ and $A_{d}$
commute, the gain matrix $\Gamma A_{d}$ symmetric bloc diagonal positive
definite. 
\begin{thm}
\label{thm: theorem RSR}Consider System (\ref{eq:real kinematics})-(\ref{eq. real dynamics}),
under assumptions in sub-section (\ref{sub:Assumtions}) and the control
law (\ref{eq:  control law quaternion}) with the auxiliary system
given by (\ref{eq:new auxiliary error dynamics}). Then, if Hypothesis
\textbf{(GEN)}. holds true, one gets that 
\begin{description}
\item [{$(1)$}] there are eight equilibrium points, given by 
\[
\Omega_{1}^{+}=(\mathbf{0}_{3n},\,\left[\begin{array}{c}
1\\
\boldsymbol{0}
\end{array}\right],\,\boldsymbol{0}),\ \Omega_{1}^{-}=(\mathbf{0}_{3n},\,\left[\begin{array}{c}
-1\\
\boldsymbol{0}
\end{array}\right],\,\boldsymbol{0}),\ \Omega_{2,3,4}^{+}=(\mathbf{0}_{3n},\,\left[\begin{array}{c}
0\\
v_{i}
\end{array}\right],\,\boldsymbol{0}),\ \Omega_{2,3,4}^{-}=(\mathbf{0}_{3n},\,\left[\begin{array}{c}
0\\
-v_{i}
\end{array}\right],\,\boldsymbol{0}),
\]
with $(v_{1},v_{2},v_{3})$ is an orthonormal basis diagonalizing
$W_{\rho}$. 
\item [{$(2)$}] All trajectories of (\ref{eq:real kinematics})-(\ref{eq. real dynamics})
converge to one of the equilibrium points defined in Item $(1)$. 
\item [{$(3)$}] Set $c:=4\lambda_{min}(W_{\rho})$, where $\lambda_{min}(W_{\rho})$
is the smallest eigenvalue of $W_{\rho}$. Then the equilibrium point
$\Omega_{1}^{+}$ is locally asymptotically stable with a domain of
attraction containing the set 
\begin{equation}
V_{c}^{+}:=\left\{ \chi\in\Upsilon\mid V(\chi)<c\;:\hbox{and}\;\bar{q}_{0}>0\right\} \label{eq:V+ attraction}
\end{equation}
and the equilibrium point $\Omega_{1}^{-}$ is locally asymptotically
stable with a domain of attraction containing the set 
\begin{equation}
V_{c}^{-}:=\left\{ \chi\in\Upsilon\mid V(\chi)<c\;\hbox{and}\;\bar{q}_{0}<0\right\} .\label{eq:V- attraction}
\end{equation}

\item [{$(4)$}] The other equilibrium points $\Omega_{2,3,4}^{\pm}$ are
hyperbolic and not stable (i.e. the eigenvalues of each of the corresponding
linear systems have non zero real part and at least one of them has
positive real part). This implies that System (\ref{eq:real kinematics})-(\ref{eq. real dynamics})
is almost globally asymptotically stable with respect to the two equilibrium
points $\Omega_{1}^{\pm}$ in the following sense: there exists an
open and dense subset $\Upsilon_{0}\subset\Upsilon$ such that, for
every initial condition $\chi_{0}\in\Upsilon_{0}$, the corresponding
trajectory converges asymptotically to either $\Omega_{1}^{+}$ or
$\Omega_{1}^{-}$. 
\end{description}
\end{thm}
\begin{IEEEproof}
Regarding Item $(1)$, one must solve the equation $f(\chi)=0$, where
$f$ is the nonlinear function describing (\ref{eq: closed loop dynamics quaternion}).
Two cases can be considered. Assume first that $\bar{q}_{0}\neq0$.
Both matrices $\bar{q}_{0}\boldsymbol{I}+S(\bar{q})$ and $\bar{q}_{0}\boldsymbol{I}-S(\bar{q})$
are non singular. Therefore from the third equation of (\ref{eq: closed loop dynamics quaternion})
$\omega=\boldsymbol{0}$ and thus $\xi=\mathbf{0}_{3n}$ from the
first equation of (\ref{eq: closed loop dynamics quaternion}). The
fourth equation of (\ref{eq: closed loop dynamics quaternion}) reduces
to $z_{\rho}=0$ and one concludes that $\bar{q}=\boldsymbol{0}$
and $\bar{q}_{0}=\pm1$ leading to two equilibrium points : $\Omega_{1}^{+}=(\mathbf{0}_{3n},\,\left[\begin{array}{c}
1\\
\boldsymbol{0}
\end{array}\right],\,\boldsymbol{0})$ and $\Omega_{1}^{-}=(\mathbf{0}_{3n},\,\left[\begin{array}{c}
-1\\
\boldsymbol{0}
\end{array}\right],\,\boldsymbol{0})$.

Assume that $\bar{q}_{0}=0$. Then $\Vert\bar{q}\Vert=1$ and according
to the third equation of (\ref{eq: closed loop dynamics quaternion}),
one gets that $\omega$ is parallel to $\bar{q}$, let say $R_{d}\omega=\mu\bar{q}$
and then $\mu$ must be equal to zero according to the second equation
of (\ref{eq: closed loop dynamics quaternion}), implying that $\omega=0$.
As in the previous case, one deduces that $\xi=\mathbf{0}_{3n}$.
The fourth equation of (\ref{eq: closed loop dynamics quaternion})
yields that $\bar{q}$ and $W_{\rho}\bar{q}$ are parallel, leading
to the six points $\Omega_{2,3,4}^{\pm}$.

We now turn to an argument for Item $(2)$. Using the facts that 
\[
\omega^{T}S(\omega)=\boldsymbol{0},\;\bar{q}^{T}W_{\rho}(\bar{q}_{0}\boldsymbol{I}+S(\bar{q}))\omega=\omega^{T}(\bar{q}_{0}\boldsymbol{I}-S(\bar{q}))W_{\rho}\bar{q},\ \omega^{T}B^{T}\Gamma\xi=\xi^{T}\Gamma B\omega,
\]
the time derivative of (\ref{eq: Lyapunov function}) in view of (\ref{eq: closed loop dynamics quaternion})
yields 
\begin{equation}
\dot{V}=-\xi^{T}\Lambda\xi\leq0,\label{eq: Lyapunov function time derivative}
\end{equation}
since $\Lambda=A_{d}^{T}\Gamma A_{d}+\Gamma A_{d}^{2}=2\Gamma A_{d}^{2}$
is symmetric positive definite. We deduce that all trajectories of
(\ref{eq: closed loop dynamics quaternion}) are defined for all times
and bounded.

Since (\ref{eq: closed loop dynamics quaternion}) is autonomous and
$V$ is radially unbounded, one can use LaSalle's invariance theorem,
cf. (\ref{eq: Lyapunov function time derivative}). Therefore every
trajectory converges to a trajectory $\gamma$ along which $\dot{V}\equiv0$.
Then $\xi$ must be identically equal to zero, implying at once that
$B_{d}\omega\equiv0$ as well. The latter assertion yields that $\omega$
must be collinear to all the $b_{i}$'s, which can be true only if
$\omega\equiv0$ since there are at least two non-collinear vectors
$b_{i}$. From the fourth equation of (\ref{eq: closed loop dynamics quaternion})
one can conclude that $z_{\rho}=\boldsymbol{0}$ leading to the conclusion
by Lemma~\ref{lem:The-solutions-of z_rho}.

We next address Item $(3)$. We provide a proof only for $\Omega_{1}^{+}$
since the other case is entirely similar. Take an initial condition
$\bar{\chi}$ in $\Omega_{1}^{+}$. Since $V$ is decreasing, the
corresponding trajectory stays in $V_{c}^{+}$ for all times and,
for every $t\geq0$, $\bar{q}(t)^{T}W_{\rho}\bar{q}(t)\leq\lambda_{min}(W_{\rho})$.
This implies that $\Vert\bar{q}(t)\Vert<1$ for every $t\geq0$ and
thus $\bar{q}_{0}(t)\neq0$ for every $t\geq0$. We deduce that $\bar{q}_{0}(t)$
keeps the same sign namely that $\bar{q}_{0}(0)$, which is positive.
Since the trajectory converges to one of the eight equilibrium points,
it must be $\Omega_{1}^{+}$ since this is the only one contained
in $V_{c}^{+}$.

We finally provide an argument for Item $(4)$. First of all notice
the equilibrium points $\Omega_{i}^{\pm}$, $i=2,3,4$, cannot be
locally asymptotically stable. Indeed let $\Omega$ be one of these
points and $U$ any open neighborhood of $\Omega$ in $\Upsilon$.
Define 
\[
V_{\Omega}^{-}:=\left\{ \chi\in\Upsilon\mid V(\chi)<V(\Omega\right\} ,
\]
and set $U^{-}:=(V_{\Omega}^{-}\cap U)$. The set $U^{-}$ is obviously
non empty since it contains points of the type $\lambda\Omega$ with
$\vert\lambda\vert<1$ close enough to $1$. Moreover, for every $\chi\in U^{-}$,
the trajectory of (\ref{eq: closed loop dynamics quaternion}) does
not converge to $\Omega$ since $V$ is non increasing.

We next prove that the linearization of (\ref{eq: closed loop dynamics quaternion})
at $\Omega$ is hyperbolic and admits an eigenvalue with positive
real part. We first perform a change of variables. If $\bar{q}_{0}=0$
then $\bar{q}=\sigma v_{\rho}$, where $\sigma=\pm1$ and $v_{\rho}$
is an eigenvector of $W_{\rho}$. Let us use the following change
of variable (cf. \cite{Chaturvedi2011,Mahony2008,Bullo2005})

\begin{equation}
X=\left[\begin{array}{c}
x_{0}\\
x
\end{array}\right]=\left[\begin{array}{c}
0\\
-\sigma v_{\rho}
\end{array}\right]\odot\left[\begin{array}{c}
\bar{q}_{0}\\
\bar{q}
\end{array}\right]=\sigma\left[\begin{array}{c}
v_{\rho}^{T}\bar{q}\\
-\bar{q}_{0}v_{\rho}-S(v_{\rho})\bar{q}
\end{array}\right].\label{eq:change_variable}
\end{equation}

From (\ref{eq:change_variable}) we have 
\begin{equation}
\left[\begin{array}{c}
\bar{q}_{0}\\
\bar{q}
\end{array}\right]=\left[\begin{array}{c}
0\\
\sigma v_{\rho}
\end{array}\right]\odot\left[\begin{array}{c}
x_{0}\\
x
\end{array}\right]=\sigma\left[\begin{array}{c}
-v_{\rho}^{T}x\\
x_{0}v_{\rho}+S(v_{\rho})x
\end{array}\right].\label{eq:change_var_inv}
\end{equation}

Rewrite (\ref{eq: closed loop dynamics quaternion}) using (\ref{eq:change_var_inv})
gives 
\begin{equation}
\begin{cases}
\begin{array}{lcl}
\dot{\xi} & = & -A_{d}\xi+B_{d}(X)\omega\\
\dot{x}_{0} & = & -\frac{1}{2}x^{T}\omega\\
\dot{x} & = & \frac{1}{2}\left(x_{0}\boldsymbol{I}+S(x)\right)\omega\\
J_{d}\dot{\omega} & = & -B_{d}^{T}(X)\Gamma A_{d}\xi-S(\omega)J_{d}\omega+2\left(x_{0}\boldsymbol{I}-S(x)\right)\left(\lambda_{\rho}I+S(v_{\rho})W_{\rho}S(v_{\rho})\right)x.
\end{array}\end{cases}\label{eq:linearized unstable system}
\end{equation}

Since the tangent space of $\mathbb{S}^{3}$ at $\left[\begin{array}{c}
1\\
\boldsymbol{0}
\end{array}\right]$ is given by the equation $x_{0}=0$, the linearization of System
(\ref{eq:linearized unstable system}) at $\Omega_{2}^{'}=(\xi,\, X,\,\omega)=(\mathbf{0}_{3n},\,\left[\begin{array}{c}
1\\
\boldsymbol{0}
\end{array}\right],\,\boldsymbol{0})$ is given by 
\begin{equation}
\dot{Z}=\mathcal{A}Z,\:\hbox{with}\:\mathcal{A}=\left[\begin{array}{ccccc}
-A_{d} &  & \boldsymbol{0} &  & H\\
\boldsymbol{0}_{3} &  & \boldsymbol{0}_{3} &  & \boldsymbol{I}_{3}/2\\
-J_{d}^{-1}H^{T}\Gamma A_{d} &  & 2J_{d}^{-1}G &  & \boldsymbol{0}_{3}
\end{array}\right],\label{eq: state space representation of unstable linearized system}
\end{equation}

where $Z=(z_{\xi}^{T},\, z_{x}^{T},\, z_{\omega}^{T})^{T}$, $z_{\xi},\, z_{x},\, z_{\omega}$are
the linearized vectors of $\xi,\, x,\,\omega$, respectively. $G=\lambda_{\rho}\boldsymbol{I}+S(v_{\rho})W_{\rho}S(v_{\rho})$
and $H=\left[\begin{array}{ccc}
H_{1}^{T} & \cdots & H_{n}^{T}\end{array}\right]^{T}$with $H_{j}=S\left(R_{d}\left(\boldsymbol{I}+2S^{2}\left(v_{\rho}\right)\right)b_{j}^{d}\right)$.
Since $\Omega$ is not locally asymptotically stable, it is enough
to show that $\mathcal{A}$ does not admit any eigenvalue with zero
real part.

Reasoning by contradiction, we thus assume that $\mathcal{A}$ has
an eigenvalue $i\, l$, $l\geq0$, with $Z^{l}=(z_{1}^{T},z_{2}^{T},z_{3}^{T})^{T}\in\mathbb{C}^{3n+6}$
a corresponding eigenvector. One gets the linear system of equations
\begin{equation}
\begin{cases}
\begin{array}{lcl}
-A_{d}z_{1}+ & Hz_{3} & =i\, lz_{1},\\
 & z_{3}/2 & =i\, lz_{2},\\
-J_{d}^{-1}H^{T}\Gamma A_{d}z_{1}+ & 2J_{d}^{-1}Gz_{2} & =i\, lz_{3}.
\end{array}\end{cases}\label{eq: linearized complexe dynamics}
\end{equation}
If $l=0$, one gets $z_{3}=z_{1}=0$ (since $A_{d}$ is positive definite)
and $J_{d}^{-1}Gz_{2}=0$. Recalling that $W_{\rho}$ is real symmetric
with distinct eigenvalues, we have that 
\[
W_{\rho}=\lambda_{\rho}v_{\rho}v_{\rho}^{T}+\lambda_{1}v_{1}v_{1}^{T}+\lambda_{2}v_{2}v_{2}^{T},
\]
where $(v_{\rho},v_{1},v_{2})$ is an orthonormal basis of $\mathbb{R}^{3}$
made of eigenvectors of $W_{\rho}$. By using the properties of $S(v_{\rho})$,
one gets that 
\[
G=\lambda_{\rho}v_{\rho}v_{\rho}^{T}+(\lambda_{\rho}-\lambda_{2})v_{1}v_{1}^{T}+(\lambda_{\rho}-\lambda_{1})v_{2}v_{2}^{T},
\]
implying that $\det(G)=\lambda_{\rho}(\lambda_{\rho}-\lambda_{1})(\lambda_{\rho}-\lambda_{2})\neq0$
and thus $z_{2}=0$. Then the eigenvector $Z$ is equal to zero, which
is impossible.

We deduce that $l>0$. One deduces that $z_{1}=(A_{d}+i\, l\boldsymbol{I}_{3n})^{-1}Hz_{3}$,
$z_{2}=-\frac{i}{2l}z_{3}$ and 
\begin{equation}
(i(J_{d}+G/l)+H^{T}\Gamma A_{d}(A_{d}+i\, l\boldsymbol{I}_{3n})^{-1}H)z_{3}=0.\label{eq:z3}
\end{equation}
Note that 
\[
H^{T}\Gamma A_{d}(A_{d}+i\, l\boldsymbol{I}_{3n})^{-1}H=\sum_{j=1}^{n}H_{j}^{T}\Lambda_{j}A_{j}(A_{j}+il\boldsymbol{I}_{3})^{-1}H_{j}.
\]
Recall that, for $1\leq j\leq n$, $A_{j}=P_{j}(\Lambda_{j})$ where
$P_{j}$ is a polynomial of degree two which is positive on $\mathbb{R}_{+}^{\ast}$.
One deduces that 
\[
H_{j}^{T}\Lambda_{j}A_{j}(A_{i}+il\boldsymbol{I}_{3})^{-1}H_{j}=\sum_{k=1}^{3}\frac{\lambda_{jk}P_{j}(\lambda_{jk})}{P_{j}(\lambda_{jk})+il}w_{jk}w_{jk}^{T},
\]
where $((H_{j}^{-1})^{T}w_{j1},(H_{j}^{-1})^{T}w_{j2},(H_{j}^{-1})^{T}w_{j3})$
is an orthonormal basis diagonalizing $\Lambda_{j}$.

Multiply Equation~\eqref{eq:z3} on the left by $(z_{3}^{\ast})^{T}$.
We get 
\[
i(z_{3}^{\ast})^{T}(lJ_{d}+G/l)z_{3}+\sum_{j=1}^{n}\sum_{k=1}^{3}\frac{\lambda_{jk}P_{j}(\lambda_{jk})(P_{j}(\lambda_{jk})-il)}{P_{j}(\lambda_{jk})^{2}+l^{2}}((z_{3}^{\ast})^{T}w_{jk})^{2}=0,
\]
where $l>0$. Since $(z_{3}^{\ast})^{T}(lJ_{d}+G/l)z_{3}$ is a real
number, we get 
\[
\sum_{j=1}^{n}\sum_{k=1}^{3}\frac{\lambda_{jk}P_{j}(\lambda_{jk})^{2}}{P_{j}(\lambda_{jk})^{2}+l^{2}}((z_{3}^{\ast})^{T}w_{jk})^{2}=0
\]
We deduce at once that $z_{3}=0$ and finally $Z=0$, which is again
a contradiction.

If $\mathcal{A}$ does not have eigenvalues with positive real part,
it would have only eigenvalues with negative real part and thus $\mathcal{A}$
would be Hurwitz, implying that (\ref{eq: closed loop dynamics quaternion})
would be locally asymptotically stable with respect to $\Omega$.
Since this is not true, we get that $\mathcal{A}$ does admit at least
one eigenvalue with positive real part. We hence proved that there
exists an unstable manifold of dimension at least one in neighborhoods
of the $\Omega_{j}^{\pm}$, $j=2,3,4$, and since all trajectories
converge to an equilibrium point we deduce that (\ref{eq: closed loop dynamics quaternion})
is almost globally asymptotically stable with respect to the two equilibrium
points $\Omega_{1}^{\pm}$. 
\end{IEEEproof}

\section{Control Gains Tuning and Simulation Results\label{sec:Simulation-results-1}}

This section provides a procedure to have ''optimal'' gains (usually
local) but approaching as near as possible to the global solution.
The effectiveness of the proposed velocity-free attitude stabilization
controller will be shown using simulation results.

We denote a state vector $\chi=(\left[\begin{array}{c}
\tilde{b}_{1}\\
\tilde{b}_{2}
\end{array}\right],\: Q,\:\omega)$, where we take two non collinear vectors $b_{1}$, $b_{2}$. For
simplicity and without loss of generality we take $R_{d}=I$, which
means that $\bar{q}=q$ and $b_{i}^{d}=r_{i}$. The matrices $\Lambda_{1}$,
$\Lambda_{2}$ are chosen diagonal such as $\Lambda_{i}=diag(\gamma_{i1},\gamma_{i2},\gamma_{i3})$
where $i=1,2$, therefore the matrices $A_{1}$, $A_{2}$ will be
$A_{i}=a_{i0}I+a_{i1}\Lambda_{i}+a_{i2}\Lambda_{i}^{2}$ where $i=1,2$.

In what follows, the following parameters are the same: the initial
angular velocity $\omega(0)=\left[0,\:0,\:0\right]^{T}$, the inertial
reference vectors $r_{1}=\left[0,\:0,\:1\right]^{T}$ and $r_{2}=\left[1,\:0,\:1\right]$,
the inertia matrix $J=diag(0.5,\:0.5,\:1)$, simulation sample time
is $0.01$s with RK4 algorithm. The notation ``TRB controller''
will be used to design the controller proposed in the paper \cite{Tayebi2013}.

\subsection{Parameters Tuning\label{sec:Parameters-Tuning-Using}}

Let us define the problem. Consider the case when we use two non collinear
inertial fixed vectors $r_{1}$, $r_{2}$ (i.e. $n=2$) and we use
the quaternion formulation of the closed loop dynamics (\ref{eq: closed loop dynamics quaternion}).
Consider now an objective function $g(\kappa)$ such that $\kappa$
is the vector of all parameters to be tuned. The problem consists
of finding $\underset{\kappa}{min}(g(\kappa))$ with the following
constraint $l(\kappa(m))\leq\kappa(m)\leq u(\kappa(m)),$ where $\kappa=[\begin{array}{cccccc}
\rho_{1} & \rho_{2} & a_{1(j-1)} & a_{2(j-1)} & \gamma_{1j} & \gamma_{2j}\end{array}]^{T}$ ($j=1,\ldots,3$, $\kappa\in(\mathbb{R}_{+}^{\ast})^{^{14}}$ is
the vector of parameters), $l(\kappa(m))$ and $u(\kappa(m))$ are
the lower and upper bounds corresponding to each parameter and $\kappa(m)$
is an element of $\kappa$.

Generally, optimization algorithms find a local optimum. It depends
on a basin of attraction of the starting point. Also, the effectiveness
of existing algorithms depends on the lower and upper limits. These
last values can be determined based on the dominant poles of the linearized
system around the stable equilibrium point.

The linearization of (\ref{eq: closed loop dynamics quaternion})
at $\Omega_{1}^{+}=(\mathbf{0}_{6},\,\left[\begin{array}{c}
1\\
\boldsymbol{0}
\end{array}\right],\,\boldsymbol{0})$ can be written as follows

\begin{equation}
\begin{cases}
\begin{array}{lcl}
\dot{z_{\xi}} & = & -Az_{\xi}+Gz_{\omega}\\
\dot{z_{q}} & = & \frac{1}{2}z_{\omega}\\
J\dot{z_{\omega}} & = & -G^{T}\Gamma Az_{\xi}-2W_{\rho}z_{q},
\end{array}\end{cases}\label{eq:linearized sys around omega+ stable}
\end{equation}
where $G=\left[\begin{array}{cc}
G_{1}^{T} & G_{2}^{T}\end{array}\right]^{T}$ with $G_{i}=S(r_{i})$, $\Gamma$ and $A$ are defined in Subsection~\ref{sub:Handling-the-lack}
and $W_{\rho}$ is defined in (\ref{eq:Wrho def}). Setting $Z=(z_{\xi}^{T},\, z_{q}^{T},\, z_{\omega}^{T})^{T}$
with $z_{\xi}\in\mathbb{R}^{6}$, $z_{q}\in\mathbb{R}^{3}$ and $z_{\omega}\in\mathbb{R}^{3}$
are the linearized vectors of $\xi,\, q,\,\omega$, respectively.
Then System (\ref{eq:linearized sys around omega+ stable}) can be
rewritten as $\dot{Z}=\mathcal{B}Z$, where 
\[
\mathcal{B}=\left[\begin{array}{ccc}
-A & \boldsymbol{0}_{3} & G\\
\boldsymbol{0}_{3} & \boldsymbol{0}_{3} & \boldsymbol{I}_{3}/2\\
-J^{-1}G^{T}\Gamma A & -2J^{-1}W_{\rho} & \boldsymbol{0}_{3}
\end{array}\right].
\]
Note that we used the fact that $z_{q_{0}}=0$. The linearization
of the closed loop dynamics is used to determine the upper and the
lower limits $u\kappa_{i}$, $l\kappa_{i}$, respectively, for each
parameter. Let's take an arbitrary initial condition in Euler angle
$[\varphi,\:\theta,\:\psi]=[30,\:10,\:45]\text{\textdegree}$ corresponding
to $Q(0)=\left[\begin{array}{c}
0.8804,\:0.2704,\:-0.02089,\:0.3891\end{array}\right]^{T}$. For an arbitrary chosen fixed $\kappa(m),\: m=3,\ldots,14$ gains
values, we start by varying $\kappa(1)$ and $\kappa(2)$. After inspecting
the zero-pole map, one can determine an upper and lower bounds for
$\kappa(1)$ and $\kappa(2)$ gains based on the placement of the
dominant pole, if it exist. Same reasoning gives the values in Table
\ref{tab:Lower-and-upper}.

\subsubsection*{Objective Functions and Optimal Control Gains Tuning}

Since there exist many possibilities to select the objective function,
we tests different objective functions derived from three well known
performance index. The first is Integral of Absolute Error (IAE),
the second is Integral of Time-weighted Absolute Error (ITAE) and
the last is Integral of Square Error, with the possibility to minimize
energy and attitude error in the same time or not, by choosing $\sigma\in[0\;1]$.
The first conclusion after several simulations is that the most appropriate
objective function for our application is the ISE function{\footnotesize{}{}
$g_{ise}(\kappa)=\int_{0}^{\infty}\left(\Vert\bar{q}\Vert^{2}+\sigma\Vert\tau\Vert^{2}\right)dt$}
with $\sigma=0.1$. Indeed, it minimizes convergence time of the quaternion
error and gives a comparable energy consumption to the ``TRB controller'',
as we will see after. Initial gain vector are chosen arbitrary as
$\kappa_{0}=[6,\:6,\:1,\:0.4,\:0.01,\:1,\:0.4,\:0.01,\:12,\:11,\:1,\:10,\:10,\:10]$.

To get an idea of the effectiveness of the optimization used methods,
we compare three functions to calculate gains optimally. The first
one uses \textit{KNITRO}, the second one is based on the use of the
Matlab \textit{fmincon} function and the third method is based on
the use of the same function as the second method with variation of
initial conditions of the parameters in a procedure called \textit{global
search} because the locality of the solution essentially depends on
the initial conditions. The best one is the third one, i.e., the \textit{global
search} method and the final value $\kappa_{final}$ with criterion
ISE is presented in Table \ref{tab:Selected gain values}. The corresponding
gain matrices are presented in Table \ref{tab:Selected gain values-1}.

\subsection{Simulation results\label{sec:Simulation-results}}

\begin{figure}[tbh]
\begin{centering}
\includegraphics[width=9cm]{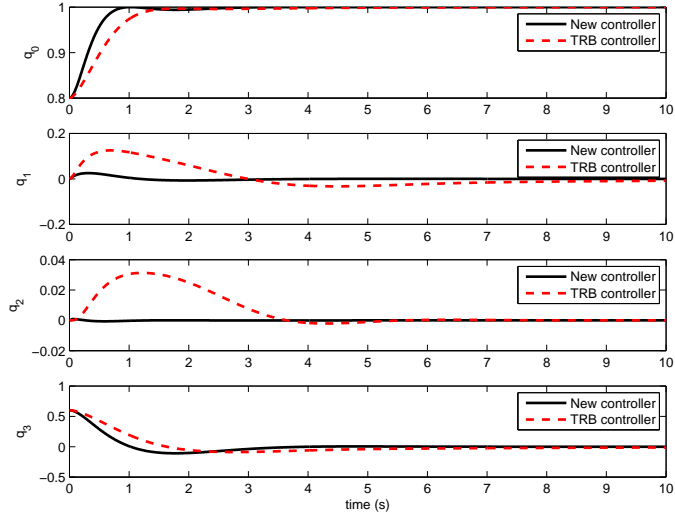} 
\par\end{centering}

\protect\caption{\label{fig:Quaternion trajectories for the two controllers}Quaternion
trajectories for the new and TRB controllers with $Q(0)=\left[\protect\begin{array}{c}
0.8,\:0,\:0,\:0.6\protect\end{array}\right]^{T}$}
\end{figure}

\begin{figure}[tbh]
\begin{centering}
\includegraphics[width=9cm]{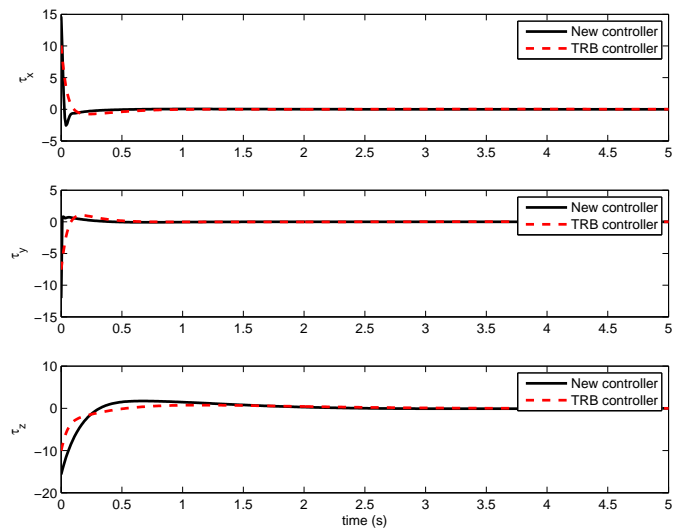} 
\par\end{centering}

\protect\caption{\label{fig:Torque for the two controllers}Applied torque for the
new and TRB controllers}
\end{figure}

Let us now show the impact of the tuned gains on the nonlinear behavior
of the new controller and the effectiveness of the proposed controller
compared with ``TRB controller''. We therefore choose the same gains
presented in \cite{Tayebi2013} for the ``TRB controller'' and the
same initial condition $Q(0)=\left[\begin{array}{c}
0.8,\:0,\:0,\:0.6\end{array}\right]^{T}$. The evolution of the unit-quaternion trajectories with respect to
time for the new and ``TRB controller'' are presented in Figure
\ref{fig:Quaternion trajectories for the two controllers}, where
the state trajectories converge asymptotically to the equilibrium
point $\Omega_{1}^{+}$. Figure \ref{fig:Torque for the two controllers}
show the torque applied in the two controllers. It is clear that the
introduction of matrix gains gives better results with a comparable
energy effort for the two controllers.

\begin{figure}[tbh]
\begin{centering}
\includegraphics[width=9cm]{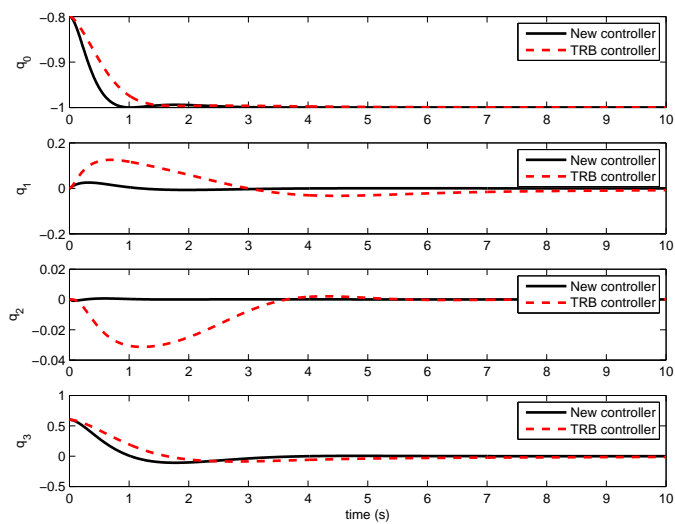} 
\par\end{centering}

\protect\caption{\label{fig:Quaternion trajectories for the two controllers q0 neg}Quaternion
trajectories for the new and TRB controllers with $Q(0)=\left[\protect\begin{array}{c}
-0.8,\:0,\:0,\:0.6\protect\end{array}\right]^{T}$}
\end{figure}

\begin{figure}[tbh]
\centering{}\includegraphics[width=9cm]{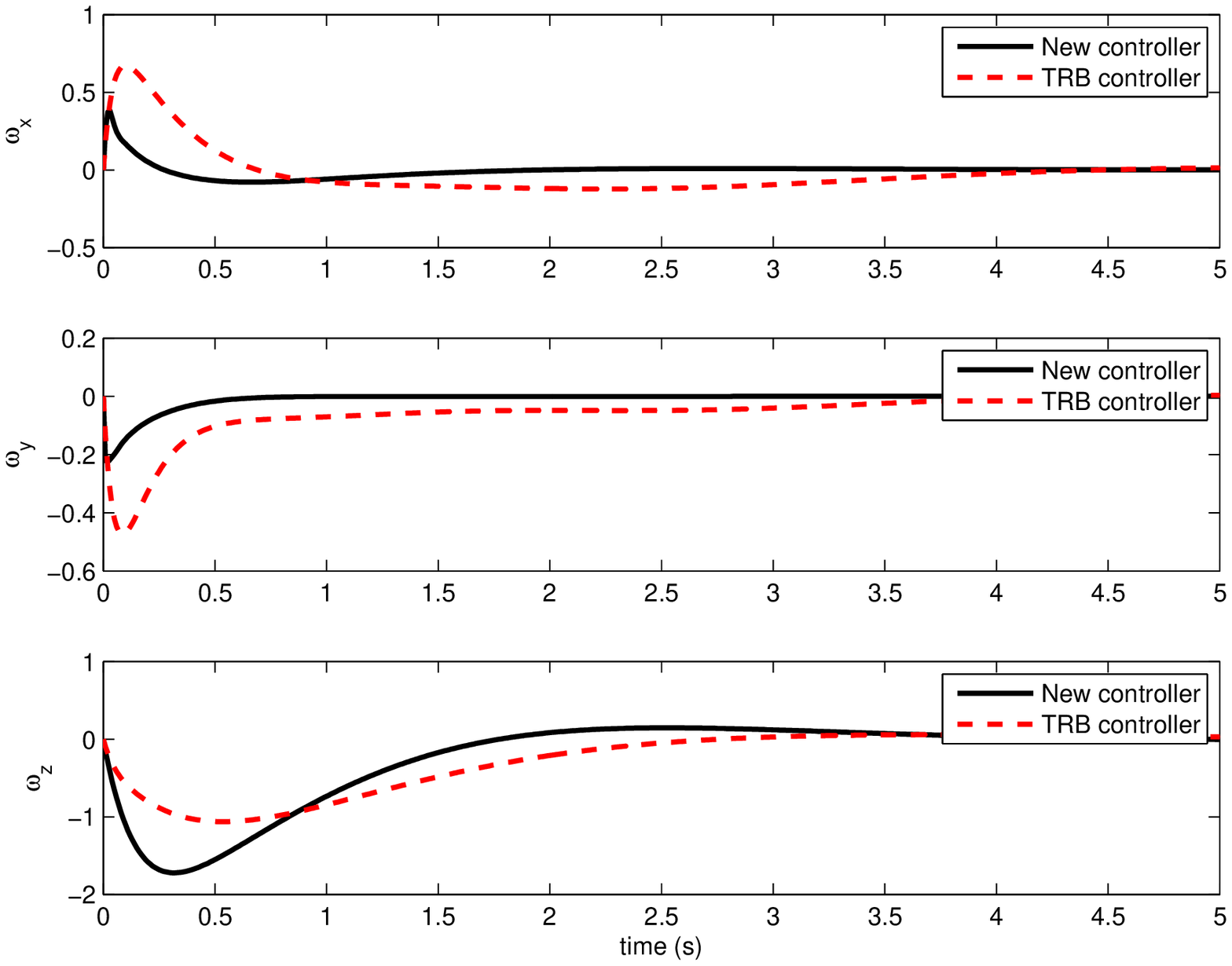}\protect\caption{\label{fig:Angular-velocity-trajectories}Angular velocity trajectories
of the new and TRB controllers. }
\end{figure}

Figure \ref{fig:Quaternion trajectories for the two controllers q0 neg}
illustrate that the proposed controller and ``TRB controller'' can
avoid the unwinding phenomenon, where the state trajectories converge
asymptotically to the equilibrium point $\Omega_{1}^{-}$ when starting
from the initial condition $Q(0)=\left[\begin{array}{c}
-0.8,\:0,\:0,\:0.6\end{array}\right]^{T}$. But, it is clear that the new controller present better performances.
Figure \ref{fig:Angular-velocity-trajectories} show the appearance
of the real angular velocity for the two controllers. 
\begin{rem}
Note that even if the initial condition is a theoretical unstable
equilibrium point, we verified by simulation that the numerical errors
push the trajectories far from this point. 
\end{rem}

\begin{rem}
The controller proposed in \cite{Thakur2014} was tested. After many
simulations, using several initial conditions, the first conclusion
is that the convergence of quaternion trajectories corresponding to
the proposed controller in the present work and ``TRB controller''
are, at least, ten time faster. The second conclusion is the fact
that the performance of the controller proposed in \cite{Thakur2014}
exhibit poor performances when only two inertial vectors are used
compared to what is presented in \cite{Thakur2014}, where results
use three vectors. 
\end{rem}

\section{Conclusions\label{sec:Conclusions}}

We have proposed an attitude stabilization controller for rigid body,
in which neither the angular velocity nor the instantaneous measurements
of the attitude are used in the feedback. This controller could be
of great help (as main or backup controllers) in applications where
prone-to-failure and expensive gyroscopes are used. When almost all
existing solutions to this problem use the instantaneous attitude
measurements, while it is well known that efficient attitude observer
use the angular velocity to obtain an accurate results, our approach
overcomes totally reconstructing the attitude. It mainly uses an auxiliary
system that can be considered as an observer of the angular velocity
using only the inertial measurements. The proposed controller doesn't
use the inertial fixed reference vectors, reduces the set of unstable
equilibria of the closed loop dynamics with respect to previous proposed
controller, provides an almost global stability of the desirable equilibrium
and avoids the \textquotedbl{}unwinding phenomenon\textquotedbl{}.
In addition, it was shown that the set of control gains leading to
a continuum of equilibria of the closed loop system is an algebraic
variety of positive co-dimension, independently on the choice of the
observed vectors. A non-linear optimal tuning method have been used
to adjust properly the controller gains. We illustrated that the introduction
of matrices gains gives a better results compared with existing work.
The performances and effectiveness of the proposed solution were illustrated
via simulation results.

%\bibliographystyle{IEEEtran}
%\bibliography{PaperARXIV}

%%%%%%%%%%%%%%%%%%%%%%%%%%%%%%%%%%%%%%%%%%%%%%%%%

\begin{table}[tbh]
\begin{centering}
\begin{tabular}{|c|c|c|}
\hline 
{\footnotesize{}gains}  & {\footnotesize{}$l(\kappa(m))$}  & {\footnotesize{}$u(\kappa(m))$}\tabularnewline
\hline 
\hline 
{\footnotesize{}$\rho_{i}(i=1,2)$}  & {\footnotesize{}0.01}  & {\footnotesize{}30}\tabularnewline
\hline 
{\footnotesize{}$a_{i0}(i=1,2)$}  & {\footnotesize{}0.0001}  & {\footnotesize{}4}\tabularnewline
\hline 
{\footnotesize{}$a_{i1}(i=1,2)$}  & {\footnotesize{}0.0001}  & {\footnotesize{}2}\tabularnewline
\hline 
{\footnotesize{}$a_{i2}(i=1,2)$}  & {\footnotesize{}0.0001}  & {\footnotesize{}0.1}\tabularnewline
\hline 
{\footnotesize{}$\gamma_{ij}(i=1,2,j=1,2,3)$}  & {\footnotesize{}0.01}  & {\footnotesize{}50}\tabularnewline
\hline 
\end{tabular}
\par\end{centering}

\protect\caption{\label{tab:Lower-and-upper}Lower and upper limits}
\end{table}

\begin{table}[tbh]
\begin{centering}
\begin{tabular}{|c|c|}
\hline 
{\footnotesize{}gains}  & {\footnotesize{}ISE $g_{ise}(\kappa)=\int_{0}^{\infty}\left(\Vert\bar{q}\Vert^{2}+0.1\Vert\tau\Vert^{2}\right)dt$}\tabularnewline
\hline 
\hline 
{\footnotesize{}$\rho_{i}(i=1,2)$}  & {\footnotesize{}{[}22.5408, 1.7736 {]}}\tabularnewline
\hline 
{\footnotesize{}$a_{1(j-1)}(j=1,2,3)$}  & {\footnotesize{}{[}4, 2, 0.1{]}}\tabularnewline
\hline 
{\footnotesize{}$a_{2(j-1)}(j=1,2,3)$}  & {\footnotesize{}{[}3.9672, 2, 0.1{]}}\tabularnewline
\hline 
{\footnotesize{}$\gamma_{1j}(j=1,2,3)$}  & {\footnotesize{}{[}50, 28.7599, 0.0971{]}}\tabularnewline
\hline 
{\footnotesize{}$\gamma_{2j}(j=1,2,3)$}  & {\footnotesize{}{[}1.8614, 1.7403, 13.9601{]}}\tabularnewline
\hline 
\end{tabular}
\par\end{centering}

\protect\caption{\label{tab:Selected gain values}Selected optimal gain values}
\end{table}

\begin{table}[tbh]
\begin{centering}
\begin{tabular}{|c|c|}
\hline 
{\footnotesize{}{}parameters}  & {\footnotesize{}{}values calculated with ISE criterion}\tabularnewline
\hline 
\hline 
{\footnotesize{}$\Lambda_{1}$}  & {\footnotesize{}diag({[}50, 28.7599, 0.0971{]})}\tabularnewline
\hline 
{\footnotesize{}$\Lambda_{2}$}  & {\footnotesize{}diag({[}1.8614, 1.7403, 13.9601{]})}\tabularnewline
\hline 
{\footnotesize{}$A_{1}$}  & {\footnotesize{}diag({[}550, 255.2727, 0.5838{]})}\tabularnewline
\hline 
{\footnotesize{}$A_{2}$}  & {\footnotesize{}diag({[}11.4541, 10.6873, 102.7916{]})}\tabularnewline
\hline 
{\footnotesize{}$W_{rho}$}  & {\footnotesize{}$\left[\begin{array}{ccc}
24.3144 & 0 & -1.7736\\
0 & 26.0881 & 0\\
-1.7736 & 0 & 1.7736
\end{array}\right]$}\tabularnewline
\hline 
{\footnotesize{}eigenvectors}  & {\footnotesize{}$v_{\rho1}=\pm\left[\begin{array}{c}
0.0780\\
0\\
0.9970
\end{array}\right]$}\tabularnewline
{\footnotesize{}of $W_{rho}$}  & {\footnotesize{}$v_{\rho2}=\pm\left[\begin{array}{c}
-0.9970\\
0\\
0.0780
\end{array}\right]$, $v_{\rho3}=\pm\left[\begin{array}{c}
0\\
-1\\
0
\end{array}\right]$}\tabularnewline
\hline 
{\footnotesize{}eigenvalues}  & {\footnotesize{}$\lambda_{\rho1}=1.6349$, $\lambda_{\rho2}=24.4531$}\tabularnewline
{\footnotesize{}of $W_{rho}$}  & {\footnotesize{} $\lambda_{\rho3}=26.0881$}\tabularnewline
\hline 
\end{tabular}
\par\end{centering}

\protect\caption{\label{tab:Selected gain values-1}Gain matrices}
\end{table}

\end{document}